\documentclass[10pt,a4paper]{amsart}
\usepackage{amsmath}
\usepackage{amssymb}
\usepackage{amsopn}
\usepackage{epsfig}
\usepackage{amsopn}
\usepackage{amsfonts}
\usepackage{latexsym}
\newtheorem{Theorem}{Theorem}
\newtheorem{Proposition}{Proposition}
\newtheorem{Lemma}{Lemma}

\newcommand{\N}{\mathbb{N}}
\newcommand{\prob}{{\rm I}\! {\rm P}}

\begin{document}

\title{Invariant densities for random $\beta$-expansions}
\author{Karma Dajani}
\address{Universiteit Utrecht,
Fac.~Wiskunde en Informatica and MRI, Budapestlaan 6, P.O. Box
80.000, 3508 TA Utrecht, the Netherlands}
\email{dajani@math.uu.nl}
\author{Martijn de Vries}
\address{Vrije Universiteit Amsterdam, Faculty of Exact Sciences, de Boelelaan 1081, 1081 HV Amsterdam, the Netherlands}
\email{mdvries@cs.vu.nl}

\maketitle
\begin{abstract}
Let $\beta >1$ be a non-integer. We consider expansions of the form $\sum_{i=1}^{\infty} d_i \beta^{-i}$, where the digits $(d_i)_{i \geq 1}$ are generated by means of a Borel map $K_{\beta}$ defined on $\{0,1\}^{\N}\times  \left[ 0, \lfloor \beta \rfloor /(\beta -1)\right]$. We show existence and uniqueness of an absolutely continuous $K_{\beta}$-invariant probability measure w.r.t. $m_p \otimes \lambda$, where $m_p$ is the Bernoulli measure on $\{0,1\}^{\N}$ with parameter $p$ $( 0 < p < 1)$ and $\lambda$ is the normalized Lebesgue measure on $[0 ,\lfloor \beta \rfloor /(\beta -1)]$. Furthermore, this measure is of the form $m_p \otimes \mu_{\beta,p}$, where $\mu_{\beta,p}$ is equivalent with $\lambda$. We establish the fact that the measure of maximal entropy and $m_p \otimes \lambda$ are mutually singular. In case the number $1$ has a finite greedy expansion with positive coefficients, the measure $m_p \otimes \mu_{\beta,p}$ is Markov. In the last section we answer a question concerning the number of universal expansions, a notion introduced in [EK].
\end{abstract}
\vskip .5cm
\section{Introduction}\label{introduction}

Let $\beta >1$ be a non-integer. In this paper we consider expansions of
numbers $x$ in $J_{\beta} := [0,\lfloor\beta\rfloor / (\beta -1)]$ of the form
$$x=\sum_{i=1}^{\infty }\frac{a_i}{\beta^i}$$
with $a_i\in\{0,1,\ldots ,\lfloor \beta \rfloor \}$, $i \in \mathbb{N}$. We shall refer to expansions of this form as $(\beta-)${\it expansions} or {\it expansions in base} $\beta$.
The largest expansion in lexicographical order of a number $x \in J_{\beta}$ is the {\it greedy expansion} of $x$; [P], [R1], [R2], and the smallest is the {\it lazy expansion} of $x$; [JS], [EJK], [DK1].
The greedy expansion is obtained by iterating the greedy transformation
$T_{\beta}: J_{\beta} \to J_{\beta}$, defined by
$$
T_{\beta}(x)\, =\, \beta x - d \quad \mbox{for} \quad x \in C(d),$$
where
$$C(j) \,=\, \left[\frac{j}{\beta},\frac{j+1}{\beta}\right),
 \quad j \in \{0, \ldots, \lfloor \beta \rfloor - 1 \},$$
and
$$C(\lfloor \beta \rfloor) \, = \, \left[\frac{\lfloor \beta \rfloor }{\beta},
\frac{\lfloor \beta \rfloor}{\beta - 1} \right].$$
The greedy expansion of $x \in J_{\beta}$ is given by $x = \sum_{i=1}^{\infty}
d_i(x)/\beta^i$, where $d_i(x)=d$ if and only if $T_{\beta}^{i-1}(x) \in C(d).$
Let $\ell: J_{\beta} \to J_{\beta}$ be given by
$$\ell(x) = \frac{\lfloor \beta \rfloor}{\beta - 1} - x,$$
then the lazy transformation $L_{\beta}: J_{\beta} \to J_{\beta}$ is defined by
$$L_{\beta}(x)= \beta x - d \quad \mbox{for} \quad
x \in \Delta(d)= \ell \left( C(\lfloor \beta \rfloor - d) \right)\, , \,
d \in \{0, \ldots, \lfloor \beta \rfloor\}.$$
The lazy expansion of $x \in J_{\beta}$ is given by
$x= \sum_{i=1}^{\infty} \tilde{d}_i(x)/\beta^i$,
where $\tilde{d}_i(x)=d$ if and only if $L_{\beta}^{i-1}(x) \in \Delta (d).$

\vskip .3cm
We denote by $\mu_{\beta}$ the extended
$T_{\beta}$-invariant {\it Parry} measure
(see [P],[G]) on $J_{\beta}$
which is absolutely continuous with respect to Lebesgue measure,
and with density
$$
h_{\beta}(x)\, =\, \left\{ \begin{array}{ll}\, \frac{1}{F(\beta )}
\sum_{n=0}^{\infty}
\frac{1}{\beta^n}\, \, 1_{[0,T_{\beta}^{n}(1))}(x)\,   & 0\leq x < 1,\\
& \\
0 & 1\leq x \leq \lfloor
\beta \rfloor /(\beta -1),
\end{array}\right.
$$
where $F(\beta )$ is the normalizing constant.
Define the {\it lazy} measure $\rho_{\beta}$ on $J_{\beta}$ by
$\rho_{\beta}=\mu_{\beta} \circ \ell^{-1}.$
It is easy to see ([DK1]) that
$\ell$ is a continuous isomorphism between $(J_{\beta},\mu_{\beta} ,T_{\beta})$ and
$(J_{\beta},\rho_{\beta} ,L_{\beta})$.

In order to produce other expansions in a dynamical way,
a new transformation $K_{\beta}$ was introduced in [DK2].
The expansions generated by iterating this map are random mixtures
of greedy and lazy expansions.
This is done by superimposing the
greedy map and the corresponding lazy map on $J_{\beta}$.
In this way one obtains
$\lfloor \beta \rfloor$ intervals on which the greedy map and
the lazy map differ. These intervals are given by
$$
S_{k} \, =\, \left[
\frac{k}{\beta},\, \frac{\lfloor \beta \rfloor }{\beta (\beta -1)}
+ \frac{k -1}{\beta} \right] \, ,\quad k =1,\ldots , \lfloor \beta
\rfloor ,
$$
which one refers to as {\em switch regions}.
On $S_{k}$, the greedy map assigns the digit $k$,
while the lazy map assigns the digit $k -1$. Outside these switch
regions both maps are identical, and hence they assign the same
digits.  Now define other expansions in base $\beta$
by randomizing the choice of the map used in the switch regions. So,
whenever $x$ belongs to a switch region, flip a coin to decide
which map will be applied to $x$, and hence which digit will be assigned.
To be more precise,  partition the interval $J_{\beta}$
into switch regions $S_{k}$ and {\em equality regions}
$E_{k}$, where
$$
E_{k}\, =\, \left( \frac{\lfloor \beta \rfloor}{\beta (\beta -1)} +
\frac{k -1}{\beta} ,\, \frac{k +1}{\beta} \right) \, ,\quad
k =1,\ldots , \lfloor \beta \rfloor -1,
$$
$$
E_0\, =\, \left[ 0,\, \frac{1}{\beta }\right) \quad {\mbox{ and
}}\quad E_{\lfloor \beta \rfloor}\, =\, \left(
\frac{\lfloor \beta \rfloor}{\beta (\beta -1)} +
\frac{\lfloor \beta \rfloor -1}{\beta} ,\,
\frac{\lfloor \beta \rfloor}{\beta -1} \right].
$$ Let
$$
S\, =\, \bigcup_{k =1}^{\lfloor \beta \rfloor} S_{k} ,
\quad {\mbox{ and }} \quad E\, =\, \bigcup_{k =0}^{\lfloor \beta \rfloor}
E_{k} ,$$
and consider $\Omega = \{ 0,1\}^{\N}$ with product
$\sigma$-algebra $\mathcal{A}.$ Let $\sigma : \Omega \to \Omega$
be the left shift, and define
$K_{\beta} : \Omega \times J_{\beta}
\to \Omega \times J_{\beta}$ by
\begin{equation*}
K_{\beta}(\omega , x)\, =\, \left\{ \begin{array}{ll}
(\omega , \beta x - k ) & x\in E_{k},\; k =0,1,\dots
,\lfloor \beta \rfloor ,\\
 & \\
(\sigma (\omega ), \beta x - k ) & x\in S_{k}\;
{\mbox{ and }}\; \omega_1=1,\, k =1,\ldots , \lfloor \beta
\rfloor ,\\
 & \\
(\sigma (\omega ), \beta x -k +1) & x\in S_{k}\;
{\mbox{ and }}\; \omega_1=0,\, k =1,\ldots , \lfloor \beta
\rfloor .
\end{array}\right.
\end{equation*}

The elements of $\Omega$ represent the coin tosses (`heads'=1 and `tails'=0)
used every time the orbit $\{K_{\beta}^n(\omega,x) : n \geq 0 \}$
hits $\Omega \times S$. Let
$$
d_1\, =\, d_1(\omega , x)\, =\, \left\{ \begin{array}{ll}
k  & {\mbox{ if }}\; x\in E_{k},\; k =0,1,\ldots , \lfloor \beta
\rfloor ,\\
 & {\mbox{ or }}\; (\omega ,x)\in \{ \omega_1=1\} \times S_{k},\;
 k =1,\ldots ,\lfloor \beta \rfloor ,\\
  & \\
k -1 &  {\mbox{ if }}\; (\omega ,x)\in \{ \omega_1=0\} \times S_{k},\;
 k =1,\ldots ,\lfloor \beta \rfloor ,
\end{array}\right.
$$
then
$$
K_{\beta}(\omega ,x)\, =\, \left\{ \begin{array}{ll}
(\omega , \beta x - d_1) & {\mbox{ if }}\; x\in E,\\
 & \\
(\sigma (\omega ), \beta x - d_1) & {\mbox{ if }}\; x\in S.
\end{array}\right.
$$
Set $d_n=d_n(\omega ,x)=d_1\left( K_{\beta}^{n-1}(\omega ,x)\right)$,
and let $\pi_2 :  \Omega \times J_{\beta}
\to J_{\beta}$ be the
canonical projection onto the second coordinate. Then
$$
\pi_2 \left( K_{\beta}^n(\omega ,x)\right) \, =\,
\beta^nx-\beta^{n-1}d_1-\cdots -\beta d_{n-1} - d_n,
$$
and rewriting yields
$$
x\, =\, \frac{d_1}{\beta} + \cdots +
\frac{d_n}{\beta^n} + \frac{\pi_2 \left( K_{\beta}^n(\omega
,x)\right)}{\beta^n} .
$$
This shows that for all $\omega \in \Omega$ and
for all $x\in J_{\beta}$  one has
that
$$
x\, =\, \sum_{i=1}^{\infty} \frac{d_i}{\beta^i}\, =\,
\sum_{i=1}^{\infty} \frac{d_i(\omega ,x)}{\beta^i}.
$$
The random procedure just described shows that with each $\omega
\in \Omega$ corresponds an algorithm that produces an expansion in
base $\beta$. Furthermore, if we identify the point $(\omega ,x)$ with
$(\omega,(d_1(\omega ,x), d_2(\omega ,x),\ldots ))$, then the action of
$K_{\beta}$ on the second coordinate corresponds to the left
shift.
\vskip .5cm

Let $<_{lex}$ and ${\leq}_{lex}$ denote the lexicographical ordering
on both $\Omega$ and $\{0,\ldots ,\lfloor \beta \rfloor \}^{\N}.$
We recall from [DdV] the following basic properties of random
$\beta$-expansions.

\begin{Theorem}\label{basic1}
Suppose $\omega ,\omega^{\prime}\in \Omega$ are such that
$\omega <_{lex} \omega^{\prime}$, then
$$(d_1(\omega ,x), d_2(\omega ,x),\ldots )\,\,
{\leq}_{\mbox{lex}}\,\, (d_1(\omega^{\prime} ,x), d_2(\omega^{\prime} ,x),
\ldots ).$$
\end{Theorem}
\medskip
\noindent

\begin{Theorem}\label{basic2}
Let $x\in J_{\beta}$ and let
$x=\sum_{i=1}^{\infty}a_i/{\beta^i}$ with $a_i\in \{0,1,\ldots,
\lfloor \beta \rfloor \}$ be an expansion of $x$ in base
$\beta$.  Then there exists an $\omega \in \Omega $ such that for all
$i \geq 1$, $a_i=d_i(\omega,x).$
\end{Theorem}

\vskip .5cm
In [DdV] it is shown that there exists a unique measure of maximal entropy
$\nu_{\beta}$ for the map $K_{\beta}$. It is the main goal of this
paper to investigate the
relationship between this measure and the measure $m_p \otimes \lambda$,
where $\lambda$ is the normalized Lebesgue measure on $J_{\beta}$ and $m_p$
is the Bernoulli measure on $\Omega = \{0,1\}^{\N}$ with parameter
$p$ $(0 <p <1)$:

$$m_p(\{\omega_1= i_1 , \ldots , \omega_n = i_n\})=
p^{\sum_{j=1}^{n} i_j} (1-p)^{n - \sum_{j=1}^{n} i_j}.$$
\vskip .3cm

In this paper, the parameter $p \in (0,1)$ is fixed but arbitrary,
unless stated otherwise.
In order to prove that the measures
$\nu_{\beta}$ and $m_p \otimes \lambda$ are mutually singular,
we introduce in the next section another
$K_{\beta}$-invariant probability measure.
This measure is a product measure $m_p \otimes \mu_{\beta,p}$ and we show
 in Section 3 that $K_{\beta}$ is ergodic w.r.t. this measure.
Furthermore, the measures $m_p \otimes \lambda$ and
$m_p \otimes \mu_{\beta,p}$ are shown to be equivalent.
These facts enable us to conclude that the measures $\nu_{\beta}$ and
$m_p \otimes \lambda$ are mutually singular. Moreover, it follows that
$m_p \otimes \mu_{\beta,p}$ is the unique absolutely continuous
$K_{\beta}$-invariant probability measure w.r.t. $m_p \otimes \lambda$.
The measure $\mu_{\beta,p}$ satisfies the important relationship
$$\mu_{\beta,p}= p \cdot \mu_{\beta,p} \circ T_{\beta}^{-1}
+ (1-p) \cdot \mu_{\beta,p} \circ L_{\beta}^{-1}.$$
In Section 4 we show that if 1 has a finite greedy expansion with
positive coefficients, then the measure $m_p \otimes \mu_{\beta,p}$
is Markov, and we determine the measure $\mu_{\beta,p}$ explicitly.
In Section 5 we discuss some open problems. As an application of some of
the results in this paper, we also show that for
$\lambda$-a.e.\ $x \in J_{\beta}$, there exist $2^{\aleph_0}$ so called
universal expansions of $x$ in base $\beta$.

\section{The skew product transformation $R_{\beta}$}
\vskip .3cm
Define the {\it skew product} transformation $R_{\beta}$ on
$\Omega \times J_{\beta}$ as follows.
$$
R_{\beta}(\omega, x)\, =\, \left\{ \begin{array}{ll}
(\sigma (\omega), T_{\beta} x) \, & {\mbox{   if  }}\,\,  \omega _1 = 1,\\

(\sigma (\omega), L_{\beta}x) \,  & {\mbox{   if  }}\,\,   \omega _1 = 0.
\end{array}\right.
$$
\noindent
On  the set $\Omega \times J_{\beta}$,
we consider the $\sigma$-algebra $\mathcal{A} \otimes \mathcal{B}$,
where $\mathcal{A}$ is the product
 $\sigma$-algebra on $\Omega$ and $\mathcal{B}$ is the Borel
$\sigma$-algebra on $J_{\beta}$. Let $\mu$ be an arbitrary probability
measure on $J_{\beta}$. It
is easy to see that $m_p \otimes \mu$ is $R_{\beta}$-invariant if and only if
$\mu = p \cdot \mu \circ T_{\beta}^{-1} + (1-p) \cdot
\mu \circ L_{\beta}^{-1}$. The following result shows that a
product measure of the form
$m_p \otimes \mu$ is $K_{\beta}$-invariant if and only if it is
$R_{\beta}$-invariant.
\begin{Lemma}\label{invariance}

$m_p \otimes \mu \circ K_{\beta}^{-1}=m_p \otimes \mu \circ R_{\beta}^{-1}
=m_p \otimes \nu$, where
$$ \nu = p \cdot \mu \circ T_{\beta}^{-1} + (1-p) \cdot
\mu \circ L_{\beta}^{-1}.$$

\end{Lemma}

\noindent
{\bf Proof}. Denote by $C$ an arbitrary cylinder in $\Omega$ and let
$[a,b]$ be an interval in $J_{\beta}$. It suffices to verify that the measures coincide on sets of the form $C \times [a,b]$, since the collection of these sets forms a generating $\pi$-system.
Furthermore, let
$[i,C]=\{\omega_1=i\} \cap \sigma^{-1}(C)$ for $i =0,1$. Note that
$E \cap T_{\beta}^{-1}[a,b] = E \cap L_{\beta}^{-1}[a,b]$, and that
\begin{eqnarray*}
K_{\beta}^{-1}(C \times [a,b]) &=& C \times (E \cap T_{\beta}^{-1}[a,b]) 
\cup \; [0,C] \times (S \cap L_{\beta}^{-1}[a,b])\\
& & \cup \; [1,C] \times (S \cap T_{\beta}^{-1}[a,b]).
\end{eqnarray*}

\noindent
Hence,
\begin{eqnarray*}
m_p \otimes \mu \circ K_{\beta}^{-1}(C \times [a,b])
&=& p \cdot m_p(C) \cdot \mu(T_{\beta}^{-1}[a,b])\\
&& + (1-p) \cdot m_p(C) \cdot \mu(L_{\beta}^{-1}[a,b])\\
&=& m_p \otimes \nu (C \times [a,b]).
\end{eqnarray*}
On the other hand,
$$
R_{\beta}^{-1}(C \times [a,b]) = [0,C] \times L_{\beta}^{-1}[a,b]
 \cup [1,C] \times T_{\beta}^{-1}[a,b],
$$
and the result follows.
\hfill $\Box$\medskip
\vskip .3cm
Let $\mathfrak{D}=\mathfrak{D}(J_{\beta},\mathcal{B},\lambda)$ denote the space of
probability density functions on $J_{\beta}$ with respect to $\lambda$.
A measurable transformation
$T : J_{\beta} \to J_{\beta}$ is called nonsingular
if $\lambda (T^{-1}B)=0$ whenever $\lambda(B)=0$.

If $\mu$ is absolutely continuous w.r.t. $\lambda$ with probability density
$f= d \mu / d \lambda$ and if $T$ is a nonsingular transformation, then
$\mu \circ T^{-1}$ is absolutely continuous w.r.t. $\lambda$
with probability density $P_{T} f$ (say). Equivalently,
the Frobenius-Perron operator $P_T : \mathfrak{D} \to
\mathfrak{D}$
is defined as a linear operator such that for $f \in \mathfrak{D}$,
$P_T f$ is the function for which
$$\int_{B} P_{T} f d \lambda = \int_{T^{-1}B} f d \lambda \quad
\mbox{for all} \quad B \in \mathcal{B}.$$
Existence and uniqueness ($\lambda$-a.e.) follow from the
Radon-Nikod\'ym Theorem. A nonsingular transformation $T : J_{\beta} \to J_{\beta}$ is said to
be a Lasota-Yorke type
map (L-Y map) if $T$ is piecewise monotone and $C^2$. Piecewise monotone
and $C^2$ means that there exists a partition $\mathcal{P}=\{[a_{i-1},a_i]
 : i=1,\ldots,k\}$, such that for each $i=1,\ldots,k$, the restriction of $T$
to $(a_{i-1},a_i)$ is monotone and extends to a $C^2$ map on $[a_{i-1},a_i]$.
For such a transformation the Frobenius-Perron operator can be computed
explicitly (see [BG, page 86]) by the formula
\begin{equation}\label{formula}
P_T f(x)= \sum_{T(y)=x} \frac{f(y)}{\vert T'(y) \vert}.
\end{equation}
If, in addition, $\vert T'(x)\vert \geq \alpha > 1$ for each $x \in (a_{i-1},a_i), \, i=1,\ldots,k$, then we say that $T$ is a piecewise expanding L-Y map.
Let $T_1 , \ldots , T_n$ be L-Y maps on $J_{\beta}$ with common partition of joint monotonicity $\mathcal{P}=\{[a_{i-1},a_i]
 : i=1,\ldots,k\}$. For
$f \in \mathfrak{D}$, define $Pf= \sum_{i=1}^{n} p_i \cdot P_{T_i}f$, where
$(p_1,\ldots,p_n)$ is a probability vector. We recall the following
important theorem, due to Pelikan; [Pel]. For more results concerning
invariant densities of L-Y maps see [LY], [LiY], [Pel].

\begin{Theorem} \label{Pelikan}

Suppose that for all $x \in J_{\beta} \setminus \{a_0, \ldots ,a_k\}$, $\sum_{i=1}^{n} \frac{p_i}{\vert T_i'(x) \vert}
\leq \gamma < 1$. Then for all $f \in \mathfrak{D}$, the limit
$$ \lim_{n \to \infty} \frac{1}{n} \sum_{j=0}^{n-1} P^{j}f=f^*$$
exists in $L_1(J_{\beta}, \lambda)$. Furthermore, $Pf^* = f^*$ and one can choose
$f^*$ to be of bounded variation.
\end{Theorem}

Since $T_{\beta}$ and $L_{\beta}$ are both piecewise expanding L-Y maps,
it follows at once from Theorem \ref{Pelikan} that for all $f \in \mathfrak{D}$,
the limit $$\lim_{n \to \infty} \frac{1}{n} \sum_{j=0}^{n-1} P^{j}f = f^*$$
exists in $L_1(J_{\beta}, \lambda)$, where
$$Pf = p \cdot P_{T_{\beta}}f + (1-p) \cdot P_{L_{\beta}}f.$$
Define for $f \in \mathfrak{D}$ the probability measure $\mu_f$ by
$$\mu_f(B)=\int_B f d \lambda \qquad [B \in \mathcal{B}].$$
Observe that $Pf = f $ if and only if
$$\mu_f = p \cdot \mu_f \circ T_{\beta}^{-1} + (1-p) \cdot
\mu_f \circ L_{\beta}^{-1},$$
{\it i.e.}, if and only if $m_p \otimes \mu_f$ is $R_{\beta}$-invariant
(cf.\ Lemma~\ref{invariance}).

Let {\bf 1} denote the constant function equal to 1 on $J_{\beta}$ and
consider the function ${\bf 1}^*$, given by
$$ {\bf 1}^* = \lim_{n \to \infty} \frac{1}{n} \sum_{j=0}^{n-1}
P^{j}{\bf 1} \quad \mbox{in } L_1(J_{\beta}, \lambda).$$

We shall assume that
the function ${\bf 1}^*$ is of bounded variation.
Note that this is possible by Theorem~\ref{Pelikan}.
It follows easily from the definition of bounded variation that the
left- and right hand limits of ${\bf 1}^*$ at every point $x \in J_{\beta}$
exist and that the function ${\bf 1}^*$ is continuous up
to countably many points. Now we modify the function ${\bf 1}^*$
in such a way that it becomes lower semicontinuous. Replace ${\bf 1}^*(x)$ at every discontinuity point $x$ in the interior of $J_{\beta}$, by setting
$${\bf 1}^*(x)= \min\{{\bf 1}^*(x^-), {\bf 1}^*(x^+)\}$$
and replace ${\bf 1}^*(x)$ by its
left- or right hand limit if $x$ is an endpoint of $J_{\beta}$.
From now on we work with this modified version of
${\bf 1}^*$ which we denote again by ${\bf 1}^*$.
The next theorem states that this function is bounded below by a
positive constant $d >0$, everywhere on $J_{\beta}$.

\begin{Theorem} \label{equivalence}

The skew product transformation $R_{\beta}$ is ergodic w.r.t. the measure
$m_p \otimes \mu_{{\bf 1}^*}$. Furthermore, the measures
$m_p \otimes \mu_{{\bf 1}^*}$ and $m_p \otimes \lambda$ are equivalent and the density ${\bf 1}^*$ is bounded below by a positive constant $d$, everywhere on $J_{\beta}$.

\end{Theorem}
\noindent
{\bf Proof}. Since $P{\bf 1}^*={\bf 1}^*$, it follows from
Lemma \ref{invariance} that the measure $m_p \otimes \mu_{{\bf 1}^*}$
is $R_{\beta}$-invariant.
It is well-known that the greedy transformation $T_{\beta}$ is ergodic
w.r.t.\ its unique absolutely continuous invariant measure, which is the Parry measure
$\mu_{\beta}$ (see Section 1).
Similarly, the lazy transformation is ergodic w.r.t.\ its
unique absolutely continuous invariant measure. This implies [Pel, Corollary 7] that the skew product transformation $R_{\beta}$ is ergodic w.r.t.\
$m_p \otimes \mu_{{\bf 1}^*}$. Since the random Frobenius-Perron operator
$P$ is integral preserving w.r.t. $\lambda$, we have that
$$\int_{J_{\beta}} {\bf 1}^* d\lambda=1.$$
In particular, there exists a point $x_0$ in the interior of $J_{\beta}$ for which
${\bf 1}^*(x_0) >0.$ By lower semicontinuity of ${\bf 1}^*$, there exist an open interval $(a,b) \subset J_{\beta}$ and a constant $c>0$ such that ${\bf 1}^*(x) > c$ for each $x \in (a,b)$. Rewriting $(\ref{formula})$ one gets

\begin{equation} \label{perfrob}
P_{T_{\beta}}f(x) = \frac{1}{\beta}\sum_{T_{\beta}y=x} f(y)
\,\, , \,\,
P_{L_{\beta}}f(x) = \frac{1}{\beta}\sum_{L_{\beta}y=x} f(y),
\end{equation}
\noindent
see also [P, Theorem 1], and thus 
$${\bf 1}^*(x) = \frac{p}{\beta} \sum_{T_{\beta}y=x} {\bf 1}^*(y) +
\frac{1-p}{\beta} \sum_{L_{\beta}y=x} {\bf 1}^*(y).$$
Hence, for $\lambda$-a.e.\ $x \in T_{\beta}(a,b)$, we have that
$${\bf 1}^*(x) > \frac{pc}{\beta}. $$
By induction, we have that for each $n$ and for $\lambda$-a.e.\
$x \in T_{\beta}^n(a,b)$,
$${\bf 1}^*(x) > \frac{p^nc}{\beta^n}. $$
It is easy to verify that there exist a number $\delta>0$ and
a positive integer $n$, such that
$$T_{\beta}^{n}(a,b) \supset [z, z+ \delta),$$
where $z$ is a discontinuity point of $T_{\beta}$.
Hence,
$$T_{\beta}^{n+1}(a,b) \supset [0, \beta \delta).$$
Moreover, there exists a positive integer $m$,
such that
$$L_{\beta}^m([0, \beta \delta)) =J_{\beta}.$$
Using the same argument as before, we conclude that for $\lambda$-a.e.\
$x \in J_{\beta}$,
$${\bf 1}^*(x) > d:= \frac{p^{n+1}(1-p)^{m} c}{\beta^{n+m+1}}.$$
Hence, the function ${\bf 1}^*$ is larger or equal than $d$ at every
continuity point of ${\bf 1}^*$. Due to our modification of ${\bf 1}^*$ at
discontinuity points, the function ${\bf 1}^*$ is everywhere larger or equal
than $d$. The equivalence of $m_p \otimes \mu_{{\bf 1^*}}$ and $m_p \otimes \lambda$ is an immediate consequence.
\hfill $\Box$\medskip
\vskip .3cm
Since any invariant probability measure absolutely continuous w.r.t.
an ergodic invariant probability measure coincides with this measure, we have
from Theorem~\ref{Pelikan} and Theorem~\ref{equivalence} that for all $f \in \mathfrak{D}$
$$\lim_{n \to \infty} \frac{1}{n} \sum_{j=0}^{n-1} P^{j}f = {\bf 1}^*
\quad \mbox{in  } \, L_1(J_{\beta},\lambda).$$

\vskip .5cm
\noindent
{\bf Remarks 1.} (1) From now on we write $\mu_{\beta,p}$
instead of $\mu_{{\bf 1}^{*}}$, since the measure depends
on both $\beta$ and $p$. It is the unique probability measure,
absolutely continuous w.r.t. $\lambda$, satisfying the relationship
\begin{equation}\label{remark}
\mu_{\beta,p} = p \cdot \mu_{\beta,p} \circ T_{\beta}^{-1}
+ (1-p) \cdot \mu_{\beta,p} \circ L_{\beta}^{-1}.
\end{equation}

(2) Recall that $\ell: J_{\beta} \to J_{\beta}$, given by $\ell(x)= \lfloor \beta \rfloor
/(\beta-1) -x $, satisfies $T_{\beta} \circ \ell = \ell \circ L_{\beta}$.
It follows from the previous remark that $\mu_{\beta,p} \circ \ell^{-1}
=\mu_{\beta, 1-p}$. In particular, we see that the invariant density
${\bf 1}^*$ is symmetric on $J_{\beta}$ if $p = 1/2$.

(3) Let $T_1, \ldots, T_n$ be piecewise expanding L-Y maps on $J_{\beta}$ and let $(p_1, \ldots , p_n)$
be a probability vector. Recently it has been shown by Boyarsky, G\'ora and
Islam (see [BGI]) that functions $f \in \mathfrak{D}$ satisfying $f=Pf=
\sum_{i=1}^n p_i \cdot P_{T_i}f$, are bounded below by a positive constant on
their support ($\lambda$-a.e.). Hence, the fact that ${\bf 1}^*$ is bounded below by a positive constant on $J_{\beta}$ can also be deduced from their result combined with the equivalence of $m_p \otimes \lambda$ and $m_p \otimes \mu_{\beta,p}$.

(4) It is well-known that the Parry measure $\mu_{\beta}$ is the unique
probability measure, absolutely continuous w.r.t.\ $\lambda$ and satisfying
equation $(\ref{remark})$ with $p=1$. Note however that $\mu_{\beta}$ and
$\lambda$ are {\it not} equivalent on $J_{\beta}$. Similarly, the lazy measure
$\rho_{\beta}$ and $\lambda$ are not equivalent. For this reason,
we restrict ourselves to values of the parameter $p$
in the open interval $(0,1)$.
\vskip .3cm

\section{Main Theorem}

It is the object of this section to show that the measure of
maximal entropy $\nu_{\beta}$ for the map $K_{\beta}$
and the measure $m_p \otimes \lambda$ are mutually singular.
\vskip .3cm
\noindent
Let $D = \{0,1,\ldots,\lfloor \beta \rfloor\}^{\mathbb{N}}$ be equipped
with the product $\sigma$-algebra $\mathcal{D}$ and let $\sigma'$ be the left
shift on $D$. Define the function $\varphi: \Omega \times J_{\beta} \to D$ by
$$\varphi (\omega,x) = (d_1(\omega,x),d_2(\omega,x),\ldots).$$
Clearly, $\varphi$ is measurable and $\varphi \circ K_{\beta} =
\sigma' \circ \varphi$. Furthermore, Theorem \ref{basic2} implies that
$\varphi$ is surjective.
\noindent
Let $$Z = \{ (\omega,x) \in \Omega \times J_{\beta} :
K_{\beta}^n(\omega,x) \in \Omega \times S \textrm{ for infinitely many }
n \geq 0 \},$$
and
$$D' = \{ (a_1,  a_2, \ldots) \in D : \sum_{i=1}^\infty \frac{a_{j + i -1}}
{\beta^i} \in S \textrm{ for infinitely many }j \geq 1 \}.$$
Observe that $K_{\beta}^{-1}(Z)=Z$, $(\sigma')^{-1}(D')=D'$ and that
the restriction
$\varphi': Z \to D'$ of the map $\varphi$ to $Z$ is a bimeasurable bijection.
Let $\prob$ denote the uniform product measure on $D$.
We recall from [DdV] that the measure $\nu_{\beta}$ defined on
$\mathcal{A} \otimes \mathcal{B}$ by
$\nu_{\beta}(A)= \prob (\varphi(Z \cap A))$ is the unique
$K_{\beta}$-invariant measure of
maximal entropy $\log (1 + \lfloor \beta \rfloor)$. It was also shown
that the projection of $\nu_{\beta}$ on the second coordinate is an
infinite convolution of Bernoulli measures (see [E1], [E2]). More precisely,
consider the purely discrete probability
measures $\{\delta_i\}_{i \geq 1}$ defined on
$J_{\beta}$ and determined by:
$$\delta_i(\{k \beta^{-i} \})=\frac{1}{\lfloor \beta \rfloor + 1}
\quad \mbox{for} \quad k=0,1, \ldots , \lfloor \beta \rfloor. $$
Let $\delta_{\beta}$ be the
corresponding infinite Bernoulli convolution,
$$\delta_{\beta}  = \lim_{n \to \infty} \delta_1 \ast \cdots \ast \delta_n ,$$
then $\nu_{\beta} \circ \pi_2^{-1} =\delta_{\beta}$.

For $\omega \in \Omega$, let $\overline{\omega}$ be given by
$$\overline{\omega}=(\overline{\omega_1},\overline{\omega_2},\ldots)=
(1-\omega_1, 1-\omega_2,\ldots).$$ Concerning the projection
$\pi_1: \Omega \times J_{\beta} \to \Omega$ of the measure
$\nu_{\beta}$ on the first coordinate, we have the following lemma.

\begin{Lemma}\label{omega}
For $n \geq 1$ and $i_1, \ldots , i_n \in \{0,1\}$, we have
$$\nu_{\beta} \circ \pi_1^{-1}(\{\omega_1 = i_1,\ldots,\omega_n = i_n\})=
\nu_{\beta} \circ \pi_1^{-1}(\{\overline{\omega_1}=i_1 ,
\ldots ,\overline{\omega_n}=i_n\}).$$
\end{Lemma}
\noindent
{\bf Proof.} Define the map $r: D \to D$ by
$$r(a_1, a_2,\ldots)=(\lfloor \beta \rfloor - a_1,\lfloor \beta \rfloor - a_2,
\ldots).$$
It follows easily by induction that for $i \geq 1$ and
$(\omega,x) \in \Omega \times J_{\beta}$,
$$d_i(\omega,x)= \lfloor \beta \rfloor - d_i(\overline{\omega},\ell(x)).$$
Hence,
$$\varphi(\omega,x)= r \circ \varphi(\overline{\omega},\ell(x)).$$
Since the map $r$ is clearly invariant w.r.t. $\prob$, the assertion follows.
\hfill $\Box$ \medskip
\vskip .3cm
In particular, it follows from Lemma \ref{omega} that
$\nu_{\beta}\circ \pi_1^{-1}(\{\omega_i = 1\})=\frac{1}{2}$, for all
$i \geq 1$. However, in general, the measure $\nu_{\beta} \circ \pi_1^{-1}$
is {\it not} the uniform Bernoulli measure on $\{0,1\}^{\mathbb{N}}$.
For instance, using the techniques in [DdV, Section 4], one easily shows that
if the greedy expansion of 1 in base $\beta$ satisfies $1= 1 / \beta +
1 / \beta^3$, then $\nu_{\beta} \circ \pi_1^{-1}$ provides a
counter example. In the case that
1 has a finite greedy expansion with positive
coefficients, it has been shown in [DdV, Theorem 8]
that $\nu_{\beta} \circ \pi_1^{-1}$ {\it is} the uniform
Bernoulli measure.
The next lemma shows that the $K_{\beta}$-invariant measures
$\nu_{\beta}$ and $m_p \otimes \mu_{\beta,p}$ are different.

\begin{Lemma}\label{different1}
$\nu_{\beta} \not= m_p \otimes \mu_{\beta,p}.$
\end{Lemma}
\noindent
{\bf Proof.} According to Theorem~\ref{equivalence}, there exists a constant $c > 0$, such that ${\bf 1^*}(x) \geq c$ for all $x \in J_{\beta}$. Choose $n \in \mathbb{N}$ such that
$\frac{1}{\beta} + \frac{1}{\beta^n} \in S_1$. Now, suppose the converse is true, {\it i.e.}, suppose that the measures $\nu_{\beta}$  and $m_p \otimes \mu_{\beta,p}$ coincide. In particular, we assume that $\nu_{\beta}$ is
a product measure and that $\delta_{\beta}=\mu_{\beta,p}$.\\
On the one hand we infer from Lemma \ref{omega} that
$$\nu_{\beta}(\{\omega_1=1\} \times J_{\beta} \,\big|\, \Omega \times
[\frac{1}{\beta},\frac{1}{\beta} + \frac{1}{\beta^n}))=
\frac{1}{2}.$$
On the other hand, since the digits $(d_i)_{i \geq 1}$ form a uniform
Bernoulli process under $\nu_{\beta}$,
\begin{eqnarray*}
&&\nu_{\beta}(\{\omega_1=1\} \times J_{\beta} \,\big|\, \Omega \times
[\frac{1}{\beta},\frac{1}{\beta} + \frac{1}{\beta^n}))\\
&=& \nu_{\beta}(\{d_1 = 1\} \vert  \Omega \times
[\frac{1}{\beta},\frac{1}{\beta} +
\frac{1}{\beta^n}))\\
&=& \frac{\nu_{\beta}(\{d_1=1 , d_2=0 ,  \ldots , d_n =0 , \sum_{i=1}^{\infty}
\frac{d_{n+i}}{\beta^i} \in [0, 1)\})}{\mu_{\beta,p}([\frac{1}{\beta},
\frac{1}{\beta} + \frac{1}{\beta^n}))}\\
&\leq& \frac{1}{c}\left( \frac{\beta}{\lfloor \beta \rfloor + 1} \right)^n
\delta_{\beta}([0,1)).
\end{eqnarray*}
Passing to the limit, we get a contradiction.
\hfill $\Box$ \medskip

\vskip .5cm

Define the map $F : \Omega \times J_{\beta} \to D$ by
$$F(\omega,x)= (d_1(\omega,x),d_1(R_{\beta}(\omega,x)),d_1(R_{\beta}^2
(\omega,x)),\ldots).$$
We have that $\sum_{i=1}^{\infty}d_1(R_{\beta}^{i-1}(\omega,x))/\beta^i = x$ for all $(\omega,x) \in \Omega \times J_{\beta}$. Moreover, the map $F$ is
surjective and
$\sigma' \circ F = F \circ R_{\beta}$. Hence $F$ is a factor map and $\sigma'$
is ergodic w.r.t. the measure $\rho = m_p \otimes \mu_{\beta,p} \circ F^{-1}$.
Note however, that the map $F$ is not injective, even if we restrict $F$ to
the set for which $R_{\beta}$ hits $\Omega \times S$ infinitely many times;
this is due to the fact that in equality regions only
one digit can be assigned.
It follows from Theorem \ref{equivalence} and Birkhoff's Ergodic Theorem
that $\rho$ is concentrated on $D'$. Therefore, the
measure $\rho'$ defined on $\mathcal{A} \otimes \mathcal{B}$ by $\rho'(A)
=\rho(\varphi(A \cap Z))$ is a $K_{\beta}$-invariant probability measure and
$K_{\beta}$ is ergodic w.r.t. $\rho'$.
\begin{Lemma} \label{equality}
$\rho'= m_p \otimes \mu_{\beta,p}.$
\end{Lemma}
\noindent
{\bf Proof.} Let
\begin{alignat}{2}
A_{00} & =\{\omega_1=0\} \times S_1 & \qquad A_{\lfloor \beta \rfloor 1}
& =\{\omega_1 = 1\} \times S_{\lfloor \beta \rfloor} \notag\\
A_{02} & = \Omega \times E_0  &  \qquad A_{\lfloor \beta \rfloor 2}
& =\Omega \times E_{\lfloor \beta \rfloor} \notag
\end{alignat}
and
\begin{eqnarray*}
A_{i0}&=&\{\omega_1=0\} \times S_{i+1}\\
A_{i1}&=&\{\omega_1=1\} \times S_i \\
A_{i2}&=&\Omega \times E_i,
\end{eqnarray*}
for $1 \leq i \leq \lfloor \beta \rfloor -1$. Note that for all $i$,
$\varphi^{-1}(\{d_1=i\})$ is the union of the sets $A_{ij}$.
It is enough to show that $\rho'= m_p \otimes \mu_{\beta,p}$
on sets of the form
$$\varphi^{-1}(\{d_1=i_1 ,\ldots, d_n = i_n\}).$$ Now,
$$\varphi^{-1}(\{d_1=i_1 ,\ldots, d_n = i_n\}) =
\bigcup_{j_1,\ldots,j_n} A_{i_1 j_1} \cap \cdots \cap K_{\beta}^{-n+1}
  A_{i_n j_n},$$
where the union is taken over all $j_1,\ldots ,j_n$ for which
$A_{i_1 j_1},\ldots , A_{i_n j_n}$ are defined.
Hence, it is enough to show that
$$\rho'(A_{i_1 j_1} \cap \cdots \cap K_{\beta}^{-n+1} A_{i_n j_n})
= m_p \otimes \mu_{\beta,p}(A_{i_1 j_1} \cap \cdots \cap
K_{\beta}^{-n+1} A_{i_n j_n}).$$
It is easy to see that the set $A_{i_1 j_1} \cap \cdots \cap
K_{\beta}^{-n+1} A_{i_n j_n}$ is a product set. Denote its projection on the
second coordinate by $V_{i_1 j_1 \ldots i_n j_n}$.\\
Define $$\mathcal{U}=\{(0,0),(\lfloor \beta \rfloor ,1)\}
\cup \{(i,j) : 1 \leq i \leq \lfloor \beta \rfloor -1, j \in \{0,1\}\}$$ and
$$\{\ell_1,\ldots, \ell_L\}=\{\ell : (i_{\ell},j_{\ell}) \in \mathcal{U}\}
\subset\{1,\ldots,n\} , \quad \ell_1 < \cdots <\ell_L.$$\\ Then,
\begin{equation} \label{een}
A_{i_1 j_1} \cap \cdots \cap K_{\beta}^{-n+1} A_{i_n j_n}
= \{\omega_1 = j_{\ell_1},\ldots,\omega_L=j_{\ell_L}\} \times
V_{i_1 j_1 \ldots i_n j_n}.
\end{equation}
Note that for all $x \in V_{i_1 j_1 \ldots i_n j_n}$,
$$F^{-1} \circ \varphi( \{\omega_1 = j_{\ell_1},\ldots,\omega_L=j_{\ell_L}\}
\times \{x\})=\{\omega_{\ell_1} = j_{\ell_1},\ldots,\omega_{\ell_L}=j_{\ell_L}
\} \times \{x\}.$$\\ Therefore, \label{2}
\begin{equation} \label{twee}
F^{-1} \circ \varphi(A_{i_1 j_1} \cap \cdots \cap K_{\beta}^{-n+1}A_{i_n j_n})=\{\omega_{\ell_1} = j_{\ell_1},\ldots,\omega_{\ell_L}=j_{\ell_L}\}
\times V_{i_1 j_1 \ldots i_n j_n}.
\end{equation}
The assertion follows immediately from (\ref{een}) and (\ref{twee}).
\hfill $\Box$ \medskip

From Theorem \ref{equivalence}, Lemma \ref{different1}, Lemma \ref{equality}
and the ergodicity of $K_{\beta}$ w.r.t. $\rho'$ and $\nu_{\beta}$, we arrive
at the following theorem.

\begin{Theorem} \label{singularity}
The measures $\nu_{\beta}$ and $m_p \otimes \lambda$
are mutually singular.
\end{Theorem}
\vskip .2cm
\noindent
{\bf Remark 2.} If $\beta \in (1,2)$ is a Pisot number, the mutual
singularity of $\nu_{\beta}$ and $m_p \otimes \lambda$ is a simple consequence
of the fact that in this case $\delta_{\beta}$ and $\lambda$ are mutually
singular (see [E1],[E2]).
\section{Finite greedy expansion of $1$ with positive coefficients, and the Markov property of the random $\beta$-expansion}

In this section we assume that the greedy expansion of $1$ in base ${\beta}$
satisfies $1=b_1/{\beta}+b_2/{\beta^2}+\cdots +b_n/{\beta^n}$ with
$b_i\geq 1$ for $i=1,\ldots,n$ and $n\geq 2$
(note that $\lfloor \beta \rfloor =b_1$). It has been shown
in [DdV] that in this case the dynamics of $K_{\beta}$ can be identified
with a subshift of finite type with an irreducible adjacency matrix.

We exhibit the measure $m_p \otimes \mu_{\beta,p}$ obtained in the previous
section explicitly. Moreover, it turns out that $K_{\beta}$ is exact w.r.t.
$m_p \otimes \mu_{\beta,p}$. The mutual singularity of
$\nu_{\beta}$ and $m_p \otimes \lambda$, {\it i.e.}, Theorem \ref{singularity},
will be derived by elementary means,
independent of the results established in the previous sections.
\noindent

The analysis of the case $\beta^2=b_1\beta +1$ needs some adjustments.
For this reason, we assume here
that $\beta^2\not=b_1\beta +1$, and refer the reader to [DdV, Remarks 6(2)]
for the appropriate modifications needed for the case $\beta^2=b_1\beta +1.$
We first recall some results obtained in [DdV] briefly, without proof.
\vskip .3cm
We begin by a proposition which plays a crucial role in finding
the Markov partition describing the dynamics of $K_{\beta}$.
\begin{Proposition}
Suppose $1$ has a finite greedy expansion of the form
$1=b_1/{\beta}+b_2/{\beta^2}+\cdots +b_n/{\beta^n}.$
 If $b_j\geq 1$ for $1\leq j \leq n$, then
\begin{itemize}
\item[(i)]  $T_{\beta}^i1= L_{\beta}^i1\in E_{b_{i+1}}$, $0 \leq i \leq n-2.$
\item[(ii)] $T_{\beta}^{n-1}1= L_{\beta}^{n-1}1=\frac{b_n}{\beta} \in S_{b_n}$
,  $T_{\beta}^n1=0$, and $L_{\beta}^n1=1.$
\item[(iii)] $T_{\beta}^i(\frac{b_1}{\beta -1}-1)=
L_{\beta}^i(\frac{b_1}{\beta -1}-1)\in E_{b_1-b_{i+1}}$, $0\leq i\leq n-2.$
\item[(iv)] $T_{\beta}^{n-1}(\frac{b_1}{\beta -1}-1)=
L_{\beta}^{n-1}(\frac{b_1}{\beta -1}-1)=
\frac{b_1}{\beta(\beta-1)}+\frac{b_1-b_n}{\beta}\in
S_{b_1-b_n+1}$,  $T_{\beta}^n(\frac{b_1}{\beta -1}-1)=
\frac{b_1}{\beta -1}-1,$ and $L_{\beta}^n(\frac{b_1}{\beta -1}-1)=
\frac{b_1}{\beta -1}$.

\end{itemize}
\end{Proposition}

\noindent

\vskip .3cm
To find the Markov chain behind the map $K_{\beta}$,
one starts by refining the partition $$
{\mathcal E}=\left\{ E_0, S_1, E_1,\ldots, S_{b_1}, E_{b_1} \right\}
$$
of $\left[0,\frac{b_1}{\beta -1}\right]$, using the orbits of $1$ and
$\frac{b_1}{\beta -1}-1$ under the transformation $T_{\beta}.$
We place the endpoints of ${\mathcal E}$ together with $T^i_{\beta}1$,
$T^i_{\beta}(\frac{b_1}{\beta -1}-1)$, $i=0,\ldots ,n-2$,
in increasing order.  We use these points to form a new partition
${\mathcal C}$ which is a refinement of ${\mathcal E}$, consisting
of intervals. We write ${\mathcal C}$ as
$${\mathcal C}=\{C_0,C_1,\ldots , C_L\}.$$
We choose $\mathcal{C}$ to satisfy the following.  For $0\leq i\leq n-2,$
\begin{itemize}
\item[-] $T_{\beta}^i1\in C_j$ if and only if $T_{\beta}^i1$ is a left
endpoint of $C_j$,
\item[-] $T_{\beta}^i(\frac{b_1}{\beta-1}-1)\in C_j$ if and only if
$T_{\beta}^i(\frac{b_1}{\beta-1}-1)$ is a right endpoint of $C_j.$
\end{itemize}

\noindent
Note that this choice is possible, since the points
$T_{\beta}^{i}1,T_{\beta}^{i}(\frac{b_1}{\beta-1}-1)$ for
$0 \leq i \leq n-2$, are all different.
From the dynamics of $K_{\beta}$
on this refinement, one reads the following properties of ${\mathcal C}$.

\vskip .3cm
\begin{itemize}
\item[{\bf p1.}] $C_{0}=\left[0,\frac{b_{1}}{\beta-1}-1 \right]$ and
$C_{L}=\left[1,\frac{b_{1}}{\beta-1}\right].$

\item[{\bf p2.}] For $i=0,1,\ldots, b_1$, $E_i$ can be written as a finite
disjoint union of the form $E_i=\cup_{j\in M_i}C_j$ with
$M_0,M_1,\ldots, M_{b_1}$  disjoint subsets of $\{0,1,\ldots, L\}$.
Further, the number of elements in $M_i$ equals the number of elements in
$M_{b_1-i}.$

\item[{\bf p3.}] For each $S_i$ there corresponds exactly one
$j\in \{0,1,\ldots ,L\}\setminus \cup_{k=0}^{b_1}M_k$ such that
$S_i=C_j.$

\item[{\bf p4.}] If $C_j\subset E_i$, then $T_{\beta}(C_j)=L_{\beta}(C_j)$
is a finite disjoint union of elements of ${\mathcal C}$, say
$T_{\beta}(C_j)=C_{i_1}\cup \dots \cup C_{i_l}.$
Since $\ell(C_j)=C_{L-j}\subset E_{b_{1}-i}$, it follows that
$T_{\beta}(C_{L-j})=C_{L-i_1}\cup \dots \cup C_{L-i_l}.$

\item[{\bf p5.}] If $C_j=S_i,$ then $T_{\beta}(C_j)=C_0$ and
$L_{\beta}(C_j)=C_L.$

\end{itemize}

\vskip .5cm

To define the underlying subshift of finite type associated with the map
$K_{\beta}$, we consider the
$(L+1)\times (L+1)$ matrix $A=(a_{i,j})$ with entries in $\{0,1\}$ defined by
\begin{equation*}
a_{i,j}\, =\, \left\{ \begin{array}{ll}
1 & {\mbox{ if }}\; i\in \cup_{k=0}^{b_{1}}M_k \mbox{ and }
\lambda(C_j \cap T_{\beta}(C_i))=\lambda(C_j),\\
0 & {\mbox{ if }}\; i\in \cup_{k=0}^{b_{1}}M_k \mbox{ and }
C_i\cap T_{\beta}^{-1}C_j= \varnothing,\\
1 & {\mbox{ if }}\; i\in \{0,\ldots ,L\}\setminus
\cup_{k=0}^{b_{1}}M_k {\mbox{ and }} j=0,L,\\
0 & {\mbox{ if }}\; i\in \{0,\ldots ,L\}\setminus
\cup_{k=0}^{b_{1}}M_k {\mbox{ and }} j\not= 0,L.
\end{array}\right.
\end{equation*}

\vskip .5cm

Let $Y$ denote the topological Markov chain (or the subshift of finite type)
determined by the matrix $A$.  That is,
$Y=\{y=(y_i)\in \{0,1,\dots ,L\}^{\N}:a_{y_i , y_{i+1}}=1 \}$.
We let $\sigma_Y$ be the left shift on $Y$.
For ease of notation, we denote by $s_1,\ldots, s_{b_1}$ the states
$j\in \{0,\ldots,L\} \setminus \cup_{k=0}^{b_{1}}M_k$ corresponding to
the switch regions $S_1,\ldots , S_{b_1}$ respectively.

For each $y\in Y$, we associate a sequence $(e_i) \in \{0,1,\ldots ,b_1\}^\N$
and a point $x\in \left[ 0,\frac{b_1}{\beta -1}\right]$ as follows. Let
\begin{equation}\label{markovdigits}
e_j\, =\, \left\{ \begin{array}{ll}
i & {\mbox{if }}\; y_j\in M_i,\\
i & {\mbox{if }}\; y_j=s_i {\mbox{ and }}y_{j+1}=0 ,\\
i-1 & {\mbox{if }}\; y_j=s_i {\mbox{ and }}y_{j+1}=L.
\end{array}\right.
\end{equation}
Now set
\begin{equation}\label{digits}
x\, =\, \sum_{j=1}^{\infty} \frac{e_j}{\beta^j}.
\end{equation}
Our aim is to define a map
$\psi :Y\rightarrow \Omega\times \left[0,\frac{b_1}{\beta -1}\right]$
that commutes the actions of $K_{\beta}$ and $\sigma_Y$.
Given $y\in Y$, equations (\ref{markovdigits}) and (\ref{digits}) describe
what the second coordinate of $\psi$ should be.  In order to be able to
associate an $\omega\in \Omega,$ one needs that
$y_i\in \{s_1,\ldots ,s_{b_1}\}$ infinitely often.
For this reason it is not possible to
define $\psi$ on  all of $Y$, but only on an invariant subset.
To be more precise, let
$$
Y^{\prime} = \left\{ y=(y_1,y_2,\ldots )\in Y :\,
y_i\in \{s_1,\ldots ,s_{b_1}\}{\mbox{ for infinitely many }}
\ i{\mbox{'s}} \right\}.
$$

Define $\psi : Y^{\prime} \to \Omega \times
\left[ 0, \frac{b_1}{\beta -1}\right]$ as follows. Let
$y=(y_1,y_2,\ldots )\in Y^{\prime}$, and define $x$ as given in
($\ref{digits}$).
To define a point $\omega \in\Omega$ corresponding to $y$, we
first locate the indices $n_i=n_i(y)$
where the realization $y$ of the Markov chain is in state $s_{r}$
for some $r \in \{ 1,\ldots, b_1 \}$. That is, let $n_1<n_2<\cdots$
be the indices such that $y_{n_i}=s_{r}$ for some $r =1,\ldots,b_1$.
Define
$$
\omega_j\, =\, \left\{ \begin{array}{ll}
1 & {\mbox{if }}\; y_{n_j+1}=0,\\
0 & {\mbox{if }}\; y_{n_j+1}=L.
\end{array}\right.
$$
Now set $\psi(y)=(\omega,x).$

\vskip .3cm
The following two lemmas reflect the fact that the dynamics of $K_{\beta}$ is essentially the
same as that of the Markov chain $Y$.

\begin{Lemma}\label{oud hulpje} Let $y\in Y^{\prime}$ be such that
$\psi(y)=(\omega ,x)$. Then,
\begin{itemize}
\item[\rm (i)]
$y_1=k$ for some $k\in \bigcup_{i=0}^{b_1}M_i \Rightarrow \; x\in C_k.$
\item[\rm (ii)]
$y_1=s_i,\, y_2=0  \; \Rightarrow \; x\in S_i$ and
$\omega_1=1$ for $i=1,\ldots,b_1$.\smallskip\
\item[\rm (iii)]
$y_1=s_i,\, y_2=L  \; \Rightarrow \; x\in S_i$ and
$\omega_1=0$ for $i=1,\ldots,b_1$.
\end{itemize}
\end{Lemma}
\vskip .3cm
\begin{Lemma}\label{isomorphism} For  $y\in Y^{\prime}$, we have
$$
\psi \circ \sigma_Y (y) \, =\, K_{\beta}\circ \psi (y).
$$
\end{Lemma}
\vskip .5cm
We now consider on $Y$ the Markov measure $Q_{\beta,p}$ with transition matrix
$P=(p_{i,j})$, given by
$$
p_{i,j}\, =\, \left\{ \begin{array}{ll}
\lambda (C_i\cap T_{\beta}^{-1}C_j)/\lambda (C_i) & {\mbox{ if }}\;
 i\in \cup_{k=0}^{b_1}M_k,\\
 & \\
p & {\mbox{ if }}\; i\in \{0, \ldots ,L\}\setminus
\cup_{k=0}^{b_1}M_k{\mbox{ and }} j=0,\\
 & \\
1-p & {\mbox{ if }}\; i\in \{0, \ldots ,L\}\setminus
\cup_{k=0}^{b_1}M_k {\mbox{ and }} j=L,\\
 & \\
0 & {\mbox{ if }}\; i\in \{0, \ldots ,L\}\setminus
\cup_{k=0}^{b_1}M_k{\mbox{ and }} j \not=0,L,
\end{array}\right.
$$ and initial distribution the corresponding stationary distribution $\pi$.
\begin{Theorem}\label{product}
$Q_{\beta,p} \circ \psi^{-1}$ is a product measure of the form
$m_p \otimes \mu$.
\end{Theorem}
\noindent
{\bf Proof}. Define the measure
$\mu$ on $\left[0, \frac{b_1}{\beta - 1}\right]$ by
$$ \mu(B)=
\sum_{j=0}^{L} \frac{\lambda (B \cap C_j)}{\lambda(C_j)}
\cdot \pi (j) \qquad [B \in \mathcal{B}].$$

\vskip .3cm
\noindent
Define the Markov partition $\mathcal{P}_0$ of
$\Omega \times \left[ 0, \frac{b_1}{\beta-1} \right]$ by
$$\mathcal{P}_0=\{\Omega \times C_j : j \in \cup_{k=0}^{b_1}M_k\}
\cup \{ \{\omega_1=i\} \times S_j : i=0,1, \, j=1,\ldots,b_1\}.$$
and let $\mathcal{P}_n =\mathcal{P}_0 \vee
K_{\beta}^{-1}\mathcal{P}_0 \vee \cdots \vee K_{\beta}^{-n} \mathcal{P}_0$.
It is straightforward to see that the inverse images of elements in
$\mathcal{P}_n$ under $\psi$ are cylinders in $Y$ and that for each element
$P \in\mathcal{P}_n$, $m_p \otimes \mu(P) =
Q_{\beta,p} \circ \psi^{-1}(P)$. It follows that $Q_{\beta,p} \circ \psi^{-1} = m_p \otimes \mu.$
\hfill $\Box$ \medskip

\vskip .3cm
\noindent
Since $P$ is an irreducible transition matrix,
$\sigma_Y$ is ergodic w.r.t. $Q_{\beta,p}$ and $\pi(i) > 0$
for all $i \in \{0,\ldots, L\}$. It follows from Lemma \ref{isomorphism} that
$K_{\beta}$ is ergodic w.r.t. $m_p \otimes \mu$. Furthermore,
it is immediately seen from the definition that $\mu$ is equivalent
with $\lambda$. Hence, the measure
$Q_{\beta,p} \circ \psi^{-1}$ is equivalent with $m_p \otimes \lambda$.
\vskip.3cm
\noindent
\begin{Proposition}\label{exact}
The map $K_{\beta}$ is exact w.r.t. $m_p \otimes \mu_{\beta,p}$.
Moreover, $\mu=\mu_{\beta,p}$.
\end{Proposition}
\noindent
{\bf Proof.} It follows from Lemma \ref{invariance} and Remarks 1(1) that
$\mu=\mu_{\beta,p}$. Since the transition matrix $P$ is also
aperiodic, $\sigma_Y$ is exact w.r.t. $Q_{\beta,p}$. It follows from Lemma
\ref{isomorphism} that $K_{\beta}$ is exact w.r.t.
$m_p \otimes \mu_{\beta,p}$.
\hfill $\Box$ \medskip
\vskip .3cm
It also follows from the above proposition that the density ${\bf 1}^*$ assumes the constant value $\pi(j)/\lambda(C_j)$ on the interval $C_j$,
$j \in \{0, \ldots,L\}$.
\vskip .3cm

\noindent
{\bf Example 1.} Let $\beta=G=\frac{1}{2}(1+ \sqrt{5})$ and let
$g=G-1=\frac{1}{2}(\sqrt{5}-1)$. Note that
$1= 1/\beta + 1/\beta^2$. In this case, we let $\mathcal{C}=
\mathcal{E}$, since $1$ and $1/(\beta-1) - 1$ are already endpoints of intervals in $\mathcal{E}$. Using the techniques in this section it is easily verified that the dynamical system $(\Omega \times J_{\beta}, \mathcal{A} \otimes \mathcal{B}, m_p \otimes \mu_{\beta,p}, K_{\beta})$ is measurably isomorphic to
the Markov chain with transition matrix $P$, given by
$$P \,= \, \left(\begin{array}{lll}
g & g^2 & 0 \\
p & 0 & 1-p \\
0 & g^2 & g \\
\end{array} \right),
$$
and stationary distribution $\pi$ determined by $\pi P = \pi$.

\vskip .3cm

It remains to prove that $Q_{\beta,p} \circ \psi^{-1}$ and $\nu_{\beta}$
are mutually singular. Since $K_{\beta}$ is ergodic w.r.t. both measures,
it suffices to show that the measures do not coincide.

\begin{Lemma} \label{different2}
$\nu_{\beta} \not= Q_{\beta,p} \circ \psi^{-1}$.
\end{Lemma}
\noindent
{\bf Proof.} We distinguish between the cases $p=1/2$ and $p \not= 1/2$.\\
Suppose $p=1/2$. On the one hand we have that for all
$i \in \{1,\ldots,\lfloor \beta \rfloor\}$
\begin{eqnarray*}
\frac{i}{\beta} + \sum_{i=2}^{\infty} \frac{d_i}{\beta^i} \in S_i
\Longleftrightarrow \sum_{i=1}^{\infty} \frac{d_{i+1}}{\beta^i} \in C_0,\\
\frac{i-1}{\beta} + \sum_{i=2}^{\infty} \frac{d_i}{\beta^i} \in S_i
\Longleftrightarrow \sum_{i=1}^{\infty} \frac{d_{i+1}}{\beta^i} \in C_L.
\end{eqnarray*}
\noindent
Using the fact that the digits $(d_i)_{i \geq 1}$ form a uniform
Bernoulli process under $\nu_{\beta}$, a simple calculation yields that
$$\nu_{\beta}(\Omega \times S) =
\frac{\lfloor \beta \rfloor}{\lfloor \beta \rfloor + 1}\cdot
\nu_{\beta}(\Omega \times C_0) +
\frac{\lfloor \beta \rfloor}{\lfloor \beta \rfloor+1}\cdot
\nu_{\beta}(\Omega \times C_L).
$$
Since $\nu_{\beta}(\Omega \times C_0)=\nu_{\beta}(\Omega \times C_L)$, it
follows that
$$\frac{\nu_{\beta}(\Omega \times S)}{\nu_{\beta} (\Omega \times C_0)}=
\frac{2 \lfloor \beta \rfloor}{\lfloor \beta \rfloor + 1}.$$
On the other hand, it follows from $\pi P = \pi$ that
$$\pi(0)=\frac{1}{\beta} \pi(0) + \frac{1}{2}(\pi(s_1)
+ \cdots +\pi(s_{b_1})).$$ Rewriting one gets
$$\frac{\pi(s_1)+ \cdots +\pi(s_{b_1})}{\pi(0)} =
\frac{Q_{\beta,p}\circ \psi^{-1}(\Omega \times S)}
{Q_{\beta,p}\circ \psi^{-1}(\Omega \times C_0)}=
\frac{2(\beta -1)}{\beta}.$$
\noindent
However, $$\frac{2(\beta -1)}{\beta} \not=
\frac{2 \lfloor \beta \rfloor}{\lfloor \beta \rfloor + 1}$$ for all
non-integer $\beta$, in particular for the $\beta$'s under consideration.\\
Suppose $p \not= 1/2$. In this case, the assertion follows from the fact
that the projection of $\nu_{\beta}$ on the first coordinate is the
uniform Bernoulli measure on $\{0,1\}^{\mathbb{N}}$ [DdV, Theorem 8].
Note that this result is applicable since 1 has a finite greedy expansion
with positive coefficients.
\hfill $\Box$ \medskip

The mutual singularity of $\nu_{\beta}$ and $m_p \otimes \lambda$
follows as before.
\noindent

\section{Open problems and final remarks}
1. We have not been able to find an explicit formula for ${\bf 1^*}$.
Recall that the Parry density $h_{\beta}=P_{T_{\beta}}h_{\beta}$ is given by
$$h_{\beta}(x) = \frac{1}{F(\beta)}
\sum_{x < T_{\beta}^n(1)} \frac{1}{\beta^n}.$$
(see Section 1). We expect that the density ${\bf 1}^*$ can be expressed
in a similar way, but now the random orbits of 1 as well as the random orbits
of the complementary point $\frac{\lfloor \beta \rfloor}{\beta -1} -1$ are
involved. Let us consider an example.
\vskip .2cm

\noindent
{\bf Example 2.} Let $p=1/2$ and $\beta = 3/2$. Note that in this case $\frac{\lfloor \beta \rfloor}{\beta -1} -1 = 1$.\\
Rewriting (\ref{perfrob}) one gets
$$P_{T_{\beta}}f(x) = \frac{1}{\beta} \sum_{i=0}^1 f(\frac{x+i}{\beta})
 \cdot 1_{[0,1)}(x) + \frac{1}{\beta} f(\frac{x+1}{\beta}) \cdot
1_{[1,2]}(x)$$ and
$$P_{L_{\beta}}f(x)= \frac{1}{\beta}f(\frac{x}{\beta})\cdot 1_{[0,1]}(x) +
\frac{1}{\beta} \sum_{i=0}^1 f(\frac{x+i}{\beta})
 \cdot 1_{(1,2]}(x).$$\\
It is easy to verify that ${\bf 1} \in \mathfrak{D}$ satisfies
$P{\bf 1}={\bf 1}$, hence ${\bf 1}^* = {\bf 1}$. It follows that
$m_{1/2} \otimes \lambda$ is $K_{3/2}$-invariant.
\vskip .3cm

2. We have not been able to give an explicit formula for
$h_{m_p \otimes \mu_{\beta,p}}(K_{\beta})$. However, in the special case that
$\beta^2 = b_1 \beta + 1$, the entropy is already calculated in [DK2]:
$$h_{m_p \otimes \mu_{\beta,p}}(K_{\beta})= \log \beta -\frac{b_1}{1+\beta^2}
\left(p \log p + (1-p) \log(1-p)\right).$$
Since in this case $\pi(s_i)=\frac{1}{1 + \beta^2}$, $i=1,\ldots , b_1$, it
follows that

$$h_{m_p \otimes \mu_{\beta,p}}(K_{\beta})= \log \beta -\mu_{\beta,p}(S)
\left(p \log p + (1-p) \log(1-p)\right).$$
One might conjecture that this formula holds in general.
\vskip .3cm

3. Fix $p \in (0,1)$. It is a direct consequence of Birkhoff's Ergodic Theorem,
Theorem \ref{equivalence} and the ergodicity of $K_{\beta}$ w.r.t.
$m_p \otimes \mu_{\beta,p}$ , that for $m_p \otimes \lambda$-a.e. $(\omega,x)
\in \Omega \times J_{\beta}$,
\begin{equation}\label{Birkhoff}
\lim_{n \to \infty}
\frac{1}{n} \sum_{i=0}^{n-1} {1}_{\Omega \times S} (K_{\beta}^i(
\omega,x)) =  \mu_{\beta,p}(S)>0.
\end{equation}
In particular, we infer from (\ref{Birkhoff}), that the set
$$G= \{ x \in J_{\beta} : x \mbox{ has a unique expansion
in base $\beta$}\}$$ has Lebesgue
measure zero, since for all $(\omega,x) \in \Omega \times G$,
$K_{\beta}^n(\omega,x) \in \Omega \times E$, for all $n \geq 0$.
Let $T_0=L_{\beta}$, $T_1 =T_{\beta}$, and let
$$N=\bigcup_{n=1}^{\infty}
\{ x \in J_{\beta} : T_{u_1} \circ \cdots \circ T_{u_n}x \in G , \mbox{ for some }
u_1,\ldots,u_n \in \{0,1\}\}.$$
Since the greedy map and the lazy map are nonsingular, $\lambda(N)=0$.
Note that $\Omega \times J_{\beta} \setminus N \subset Z$ and that for
$x \in J_{\beta} \setminus N$, different elements of $\Omega$ give rise to
different expansions of $x$ in base $\beta$. We conclude that for
$\lambda$-a.e. $x \in J_{\beta}$,
there exist $2^{\aleph_0}$ expansions of $x$ in base $\beta$. For a more elementary proof
of this fact in case $\beta \in (1,2)$, we refer to [S1].
\vskip .3cm

4. Erd\H{o}s and Komornik introduced in [EK] the notion of universal
expansions. They called an expansion $(d_1,d_2,\ldots)$ in base $\beta$ of some $x \in J_{\beta}$
universal if for each (finite) block $b_1 \ldots b_n$ consisting of digits in the set $\{0,\ldots,  \lfloor  \beta \rfloor\}$, there exists an index
$k \geq 1$, such that $d_k \ldots d_{k+n-1}=b_1 \ldots b_n$. They proved that there exists a number $\beta_0 \in (1,2)$, such that for each $\beta \in (1, \beta_0)$, {\it every} $x \in (0, 1/(\beta-1))$ has a universal expansion in base $\beta$.
Subsequently, Sidorov proved in [S2] that for a given  $\beta \in (1,2)$ and for $\lambda$-a.e.\ $x \in J_{\beta}$, there
exists a universal expansion of $x$ in base $\beta$. We now strengthen his result and the conclusion of the preceding remark by the following theorem.
\begin{Theorem}\label{universal}
For any non-integer $\beta >1$, and for $\lambda$-a.e.\ $x \in J_{\beta}$, there exist $2^{\aleph_0}$ universal expansions of $x$ in base $\beta$.
\end{Theorem}

In order to prove Theorem~\ref{universal} we need the following lemma.

\begin{Lemma} \label{groterdannul} Let $\beta > 1$ be a non-integer and let
$p \in (0,1)$. Then, for $n \geq 1$ and
$i_1, \ldots, i_n \in \{0,\ldots,  \lfloor \beta \rfloor\}$, we have that
$$m_p \otimes \mu_{\beta,p}(\{d_1=i_1 , \ldots, d_n= i_n\})>0.$$
\end{Lemma}
\noindent
{\bf Proof.} By Theorem~\ref{equivalence}, it suffices to show that
$$m_p \otimes \lambda (\{d_1=i_1 , \ldots, d_n= i_n\})>0.$$
It is easy to verify that there exists a sequence $(j_1,j_2, \ldots)\in D$,
starting with $i_1 \ldots i_n$, such that the numbers $x_1, \ldots, x_n$,
given by
$$x_r = \sum_{i=1}^{\infty} \frac{j_{i+r-1}}{\beta^i}, \quad r=1,\ldots,n,$$
are elements of $J_{\beta} \setminus \partial(S)$, where $\partial(S)$
denotes the boundary of $S$.
Consider for $m \geq 1$, the set
$$I_m= \left[\sum_{i=1}^{n+m} \frac{j_i}{\beta^i}, \sum_{i=1}^{n+m}
\frac{j_i}{\beta^i} + \sum_{i=n+m+1}^{\infty}
\frac{\lfloor \beta \rfloor}{\beta^i} \right].$$
Let $y \in I_m$ and let $(a_1,a_2,\ldots)$ be an expansion $y$,
starting with $j_1 \ldots j_{n+m}$. Define
$$y_r = \sum_{i=1}^{\infty} \frac{a_{i+r-1}}{\beta^i}, \quad r=1,\ldots,n.$$  Choose $m$ large enough, so that for each $r=1,\ldots,n$, $x_r$ and $y_r$
are elements of the same equal or switch region, regardless of the values of
the digits $a_{\ell},\, \ell >n+m$, and hence regardless of the
chosen element $y \in I_m$. Note that this is possible because
$x_r \notin \partial(S)$ for $r=1, \ldots,n$.
Denote the set of indices $r\in \{1,\ldots,n\}$ for which
$x_r \in S$ by $\{\ell_1,\ldots , \ell_L\}$. Then, for suitably chosen
$u_1,\ldots,u_L \in \{0,1\}$, we have that
$$\{\omega_1=u_1, \ldots, \omega_L=u_L\} \times I_m
\subset \{d_1=i_1,\ldots,d_n=i_n\}$$
and the conclusion follows.
\hfill $\Box$ \medskip

\noindent
{\bf Proof of Theorem~\ref{universal}.} Fix $p \in (0,1)$ and let
$b_1 \ldots b_n$ be an arbitrary block. Using Birkhoff's Ergodic Theorem,
Theorem~\ref{equivalence}, Lemma~\ref{groterdannul} and the ergodicity of
$K_{\beta}$ w.r.t. $m_p \otimes \mu_{\beta,p}$, we may conclude
that for $m_p \otimes \lambda$-a.e.\ $(\omega,x)
\in \Omega \times J_{\beta}$, the block $b_1\ldots b_n$ occurs in
\begin{equation}\label{occur}
(d_1(\omega,x),d_2(\omega,x),\ldots)
\end{equation}
with positive limiting frequency
$m_p \otimes \mu_{\beta,p} (\{d_1=b_1, \ldots , d_n = b_n \}).$
In particular, we have that for $m_p \otimes \lambda$-a.e.\
$(\omega,x) \in \Omega \times J_{\beta}$,  the block $b_1 \ldots b_n$
occurs in $(\ref{occur})$. Since there are only countably many blocks,
we have that for
$m_p \otimes \lambda$-a.e.\ $(\omega,x) \in \Omega \times J_{\beta}$,
the expansion $(\ref{occur})$ is universal in base $\beta$.
An application of Fubini's Theorem yields
that there exists a Borel set
$B \subset J_{\beta} \setminus N$ of full Lebesgue measure and there exist
sets
$A_x \in \mathcal{A}$ with $m_p(A_x)=1 \, (x \in B)$, such that for all
$x \in B$ and $(\omega,x) \in A_x \times \{x\}$, the expansion
$(\ref{occur})$ is universal in base $\beta$. Since the sets $A_x$ have
necessarily the cardinality of the continuum and since different elements
of $\Omega$ give rise to different
expansions of $x$ in base $\beta$ for any $x \in J_{\beta} \setminus N$,
the assertion follows.
\hfill $\Box$ \medskip
\vskip .3cm

5. An expansion $(a_1,a_2,\ldots)$ in base $\beta$
of some number $x \in J_{\beta}$ is called {\it normal} if each block
$i_1 \ldots i_n$ with digits in $\{0, \ldots , \lfloor \beta \rfloor\}$
occurs in $(a_1,a_2, \ldots)$ with limiting frequency
$(\lfloor \beta \rfloor + 1)^{-n}.$
Note that a normal expansion is in particular universal.

Fix $p \in (0,1)$. Since $\nu_{\beta} \not= m_p \otimes \mu_{\beta,p}$
and since both measures
$\nu_{\beta}$ and $m_p \otimes \mu_{\beta,p}$ are concentrated on $Z$,
there exists a block $i_1 \ldots i_n$ such that
$$m_p \otimes \mu_{\beta,p} (\{d_1 = i_1, \ldots , d_n = i_n \}) \not=
(\lfloor \beta \rfloor +1)^{-n}.$$
Hence, for $m_p \otimes \lambda$-a.e.\
$(\omega,x) \in \Omega \times J_{\beta}$,
the expansion $(\ref{occur})$ is universal but {\it not} normal.
On the other hand,
Sidorov proved in [S2], that there exists a Borel set $V \subset (1,2)$ of
full Lebesgue measure, such that for each $\beta \in V$ and for
$\lambda$-a.e.\ $x \in J_{\beta}$,
there exists a normal expansion of $x$ in base $\beta$.

\vskip .3cm

{\it Acknowledgements}. We are very grateful to the referee of [DdV] for
suggesting the problem discussed in this paper. We also thank Pawel G\'ora, Pierre Liardet and Boris Solomyak for very fruitful discussions concerning the proof of Theorem~\ref{equivalence}. Finally, we would like to thank the referee of this paper for many helpful suggestions concerning the presentation of this paper.

\end{document}